\documentclass[11pt]{article}
\usepackage{fullpage}

\usepackage{amsfonts}
\usepackage{amsmath}
\usepackage{amssymb}
\usepackage{amsthm}
\usepackage{float}
\usepackage{mathrsfs}
\usepackage{graphicx}
\usepackage{graphics}
\usepackage{subfigure}
\usepackage[colorlinks,linkcolor=black,anchorcolor=black,citecolor=black,hyperindex=black,CJKbookmarks=black]{hyperref}

\begin{document}

\def \N {\mathbb{N}}
\def \Z {\mathbb{Z}}
\def \Q {\mathbb{Q}}
\def \R {\mathbb{R}}
\def \C {\mathbb{C}}
\def \Cbar {\overline{\mathbb{C}}}
\def \D {\mathbb{D}}
\def \Dr {\mathbb{D}_{r}}
\def \L {\mathcal{L}}
\def \M {\mathcal{M}}
\def \T {\mathcal{T}}
\def \O {\mathcal{O}}
\def \F {\mathcal{F}}
\def \zbar {\overline{z}}
\def \eps {\varepsilon}
\def \epsg {\varepsilon > 0}
\newtheorem{theorem}{Theorem}
\newtheorem{lemma}[theorem]{Lemma}
\newtheorem{proposition}{Proposition}[section]

\theoremstyle{definition}
\newtheorem{definition}[theorem]{Definition}
\newtheorem{example}[theorem]{Example}
\newtheorem{xca}[theorem]{Exercise}

\theoremstyle{remark}
\newtheorem{remark}[theorem]{Remark}
\newtheorem{corollary}[theorem]{Corollary}

%
%
%
%

\title{A note on the run length function for intermittent maps}

\author{{\bf\large Hongfei Cui}\hspace{2mm}\vspace{1mm}\\
{\it\small Wuhan Institute of Physics and Mathematics, Chinese Academy of Sciences}, {\it\small Wuhan  430071,  China}\\
{\it\small E-mail:  cui@wipm.ac.cn}\vspace{1mm}\\
{\bf\large Lulu Fang}\\
{\it\small School of Mathematics},\\ {\it\small Sun Yat-sen university}, {\it\small Guangzhou  510275,  China}\\
{\it\small E-mail:  fanglulu1230@163.com}\vspace{1mm}\\
{\bf\large Yiwei Zhang}\footnote{Corresponding author, and every author contributes equally.}\hspace{2mm}\vspace{1mm}\\
{\it\small School of Mathematics and Statistics, Center for Mathematical Sciences},\\
{\it\small Hubei Key Laboratory of Engineering Modeling and Scientific Computing},\\
{\it\small Huazhong University of Sciences and Technology},
{\it\small Wuhan 430074,  China}\\
{\it\small E-mail: yiweizhang@hust.edu.cn}\vspace{1mm}}
\date{\today }
\maketitle


\begin{abstract}
We study the run length function for intermittent maps. In particular, we show that the longest consecutive zero digits (resp. one digits) having a time window of polynomial (resp. logarithmic) length. Our proof is relatively elementary in the sense that it only relies on the classical Borel-Cantelli lemma and the polynomial decay of intermittent maps. Our results are compensational to the Erd\H{o}s-R\'{e}nyi law obtained by Denker and Nicol in \cite{dennic13}.
\end{abstract}

\bigskip

\section{Introduction}

Consider a piecewise monotone interval map $T:[0,1)\to[0,1)$ with a countable or finite partition $\{I_{j}\}$, preserving a probability measure $\nu$. Given an $x\in[0,1)$ and $k\in\N$, set $\varepsilon_{k}(x)=j$ if $T^{k-1}(x)\in I_{j}$ for some $j$.
The \emph{run length function} for this system is defined as the maximal length of consecutive $j$ digits in the sequence $(\varepsilon_{1}(x),\cdots, \varepsilon_{n}(x))$. Namely,
$$
r_n(x,j)=\max\{k\geq 0:\varepsilon_{i+1}(x)=\cdots=\varepsilon_{i+k}(x)=j~~~\mbox{for some}~~0\leq i\leq n-k\}.
$$

In this note, we concern about the run length function for a class of intermittent maps, and show the existence of an appropriate scaling length $k(n)$, to quantify the asymptomatic behavior of $r_{n}(x,j)$ for $\nu$-typical $x$.

The studies on the run length functions are motivated from both theoretical and practical aspects. First, it has been revealed that run length function of piecewise linear interval maps has several probabilistic interpretations, which have great applications in DNA string machine \cite{AGW90}, reliability theory and non-parametric statistics (see e.g., \cite{Bal01,Bar92,Bat48,Mus00,Nov01,Rev90} and the references therein). For example, when specializing $T(x) = 2x({\rm mod}~1)$ and $\nu$ being the Lebesgue measure on $[0,1)$, the unit interval endows with a finite partition $\{I_0,I_1\}$, where $I_0 =[0,1/2)$ and $I_1=[1/2,1)$. The run length function for such a system corresponds to the longest length of consecutive terms of ``heads" in a mathematical experiment of coin tossing (originally back to de Moivre Problem LXXIV in the year of 1738 \cite{deMo1738}). Indeed, Erd\H{o}s and R\'{e}nyi \cite{Eroren70} obtained that for Lebesgue almost surely $x \in [0,1)$, both $r_n(x,0)=r_n(x, 1)$ increase to infinity with the logarithmic speed $\log_2 n$ as $n$ goes to infinity.

The above discussions of run length function has been generalized into $\beta$-transformations ($\beta>1$)\cite{Tonyuzhao}, that is, $T(x) = \beta x({\rm mod}~1)$ and $\nu$ is the Parry measure on $[0,1)$ (see \cite{Par60}). In this case, the unit interval endows with a finite partition $I_0 = [0,1/\beta),I_1 = [1/\beta,2/\beta), \cdots, I_{\lceil\beta\rceil-1} = [(\lceil\beta\rceil-1)/\beta,1)$, where $\lceil\cdot\rceil$ denotes the smallest integer no less than some number. As an important by-product, the length of cylinder can be partially estimated in terms of run length function. This fact is critical in understanding many dynamical problems of $\beta$-transformations (see \cite {FW12, FWL16}). Unfortunately, comparing to $\beta$-transformations, relatively limited work has been done in studying the run length function for interval maps with non-linear components. This forms the first goal of the note.

\medskip

Besides, the run length function also associates with the so-called Erd\H{o}s-R\'{e}nyi law. Given an observation $\varphi$ in $L^{1}(\nu)$, denote by $S_{n}(\varphi):=\sum_{i=1}^{n-1}\varphi\circ T^{i}$ and
\begin{align}\label{equ_edos-renyi}
 \theta(n,K(n),\varphi)&:=\max_{0\leq i\leq n-K(n)}\{S_{i+K(n)}(\varphi)-S_{i}(\varphi)\}\nonumber\\
 &=\max\{S_{K(n)}(\varphi)\circ T^{i}:0\leq i\leq n-K(n)\},~~~~~~\forall n\in\N,
\end{align}
i.e., the maximal Birkhoff sum gaining over a window of length $K(n)$ up to time $n$. The \emph{Erd\H{o}s-R\'{e}nyi law} concerns about the existence of an appropriate scaling length $K(n)$, such that the asymptotic behavior of the fraction $\frac{\theta(n,K(n),\varphi)}{K(n)}$ has a non-degenerate limit for $\nu$-typical $x$. Inspired by the work of Erd\H{o}s and R\'{e}nyi for i.i.d random valuables \cite{Eroren70}, it is believed that under certain hyperbolicity of $T$ and regularity of $\varphi$ hypothesis, the behavior of $\theta(n,K(n),\varphi)$ asymptotically admit a dichotomy on exhibiting either ordinary law of large numbers or asymptotic law of long length function, subject to the scaling length of $K(n)$.

Initialized by the work of Grigull \cite{Gri93} for hyperbolic rational maps, there are a number of results on the studies of Erd\H{o}s-R\'{e}nyi law(as well as answering the above dichotomy), particularly for the dynamical systems satisfies the large deviation principle. For examples, see the work of Chazottes and Collet \cite{Chacol01} for uniformly expanding interval maps, and the work of Denker and Kabluchko \cite{denzab07} for Gibbs-Markov dynamics.

Later, Denker and Nicol \cite{dennic13} provide several extensions of Erd\H{o}s-R\'{e}nyi law to the non-uniformly expanding dynamical systems, including logistic-like maps and intermittent maps. Unfortunately, if the dynamical systems are absent of large derivation principle, only partial results can be obtained, and the regularity of the observation function $\varphi$ has to be assumed Lipschitz/H\"{o}lder continuous (rather than integrable). Under this framework, the characteristic functions are inapplicable to Denker and Nicol's result, and thus no results on the run length function have been known. Moreover, their proofs (e.g., \cite[Theo 4.1]{dennic13} for the intermittent maps) heavily rely on the previous dynamical Borel-Cantelli lemma results obtained Gouez\"{e}l \cite[Theo 1.1]{lesGou07}, and large deviation results obtained by Pollicott and Sharp \cite[Theo 4]{Polsha09}. Our second goal is aiming to fulfill this gap of the run length function for the intermittent maps with a proof in a more elementary way (namely, without using the two technical machineries above). This will in our belief provide a more refined description in this research direction.

\medskip

Let us now state our results more mathematically. Consider $0<\alpha<1$ and recall the \emph{intermittent maps} $T_{\alpha}: [0,1) \to [0,1)$ is defined as
\[
T_{\alpha}(x)=
\begin{cases}
x(1+2^\alpha x^\alpha),  &\text{if $0\leq x <1/2$};\\
2x-1,  & \text{if $1/2 \leq x <1$}.
\end{cases}
\]
It is well known that $T:=T_{\alpha}$ has a finite absolutely continuous invariant probability measure $\mu$.  Let $I_0 = [0,1/2)$ and $I_1 = [1/2,1)$, and we are accordingly concerned with the length of the longest consecutive zero digits in $\left( \varepsilon_1(x), \varepsilon_2(x), \cdots, \varepsilon_n(x) \right)$, and the length of the longest consecutive one digits for $\mu$-almost all $x \in [0, 1)$. Set
\[
r^\alpha_n(x)= \max\Big\{k \geq 0: \varepsilon_{i+1}(x)=\cdots=\varepsilon_{i+k}(x)=0\ \text{for some}\ 0 \leq i\leq n-k\Big\},
\]
and
\[
R^\alpha_n(x)= \max\Big\{k \geq 0: \varepsilon_{i+1}(x)=\cdots=\varepsilon_{i+k}(x)=1\ \text{for some}\ 0 \leq i\leq n-k\Big\}.
\]
With this conventions, our main theorem is stated as follows.
\begin{theorem}
For $\mu$-almost all $x \in [0,1)$,
\begin{equation}\label{equ_polynimal}
\lim_{n \to \infty} \frac{\log r^\alpha_n(x)}{\log n^{\alpha}} =1,
  \end{equation}
  and
 \begin{equation}\label{equ_exponentially}
\lim_{n \to \infty} \frac{R^\alpha_n(x)}{\log_2n} =1.
  \end{equation}
\end{theorem}

Recall that $S_n(\chi_{I_0}(x))$ is the number of positive numbers  $0 \leq i \leq n-1$ such that $T^i(x) \in I_0$, i.e.,
$$
S_n(\chi_{I_0}(x)) = \sum_{i=0}^{n-1} \chi_{I_0}\circ T^i(x),
$$
where $\chi_{I_0}$ is the indicator function of the interval $I_0$. One can easily obtain the following corollary, which says that for a particular observable the limit of the maximal average over a time window of the length $n^{\alpha_1}$ (with $\alpha_1 <\alpha$) is equal to $1$ $\mu$-almost surely.

\begin{corollary}
If $\alpha_1<\alpha$, then for $\mu$-almost $x \in [0, 1)$,
$$
\lim_{n\to \infty} \frac{\max \{ S_{n^{\alpha_1}} (\chi_{I_0} \circ T^i(x)) :  \ 0\leq i\leq n-n^{\alpha_1} \} }  {n^{\alpha_1}}= 1.
$$
\end{corollary}
Similarly, we can obtain
\begin{corollary}
Let $S_n(\chi_{I_1}(x)) = \sum_{i=1}^{n} \chi_{I_1}\circ T^i(x)$.
Then for any integer sequence  $k(n)$ that satisfying $\limsup_{n\to \infty} \frac{k(n)}{\log_2 n} \leq 1$, and for $\mu$-almost $x \in [0, 1)$, we have
$$
\lim_{n\to \infty} \frac{\max \{ S_{k(n)} (\chi_{I_1} \circ T^i(x)) :  \ 0\leq i\leq n-k(n) \} }  {k(n)}= 1.
$$
\end{corollary}

The theorem indicates that polynomial length consecutive digits is typical for the run length function of intermittent maps, and different order of the run length function can coexist (differing from the phenomenon in $\beta$-transformations). The result implies the \cite[Theo 4.1(b)]{dennic13} for $\phi =\chi_{I_0}$ a indicator observable of $I_0$.

Our proof relies heavily  on the classical Borel-Cantelli lemma and the polynomial decay of intermittent maps. We obtain the upper limit by using the equivalence between the absolutely continuous invariant measure $\mu$ and Lebesgue measure outside a small one-side neighborhood of $0$. To get a result of other direction, we break the first $n$ iterations into disjoint blocks with a small scale, on which the digits are not  all $0$ or $1$. Then using the correlations between different blocks, and combining the measure estimation of the consecutive $0$ or $1$ digits, we obtain the lower limit of the run length function.

\section{Proof of theorem}

Note that $T_{\alpha}$ has a finite absolutely continuous invariant probability measure $\mu$
with a density function $h(x) = d\mu/dx$ satisfying
\[
\lim_{x \to 0} x^\alpha h(x) =c
\]
for some constant $c$ (see Theorem A of Hu \cite{lesHu04}). As a consequence, we obtain there exists a constant $1<C<\infty$ such that
\begin{equation}\label{mu}
C^{-1}x^{1-\alpha} \leq \mu([0,x)) \leq Cx^{1-\alpha}
\end{equation}
for sufficiently small $x \in [0,1)$. Moreover, $C \leq h(x) \leq C x^{1-\alpha}$ for $x \in [a,1)$ with $a>0$, then the measure $\mu$ and Lebesgue measure are uniform equivalent on every interval $[a, 1)$, see \cite{Lsv99, lesYoung}. It follows that there exists a constant $\tilde{C}< \infty$ such that for any interval $A \subset [0,1)$, we have
\begin{equation}\label{mu0}
\frac{1}{\tilde{C}} |[1/2, 1) \cap A|\leq \mu( [1/2, 1) \cap A) \leq \tilde{C}|A|.
\end{equation}

We now turn to the proof of the theorem.
\begin{proof}
The proof is divided into three  parts:
\\
{\bf Part I.} We will prove that for $\mu$-almost all $x \in [0,1)$,
\[
\limsup_{n \to \infty} \frac{\log r^\alpha_n(x)}{\log n} \leq \alpha.
\]
Let $0<\varepsilon<1-\alpha$. It suffices to  show that $\mu \{x\in [0,1): r^\alpha_n(x) \geq n^{\alpha+\varepsilon}\ i.o.\} =0$. In fact, if $r^\alpha_n(x) \geq n^{\alpha+\varepsilon}$ holds for some $x \in [0,1)$ and infinitely many $n \in \mathbb{N}$, then we obtain two cases that either $\varepsilon_{i}(x)=0$ for all $i$, $1\leq i \leq \lceil n^{\alpha+\varepsilon}\rceil$ and infinitely many $n$,   or there exists $1\leq i\leq n-\lceil n^{\alpha+\varepsilon}\rceil$ such that
\[
\varepsilon_{i}(x)=1 \ \ \text{and} \  \varepsilon_{i+1}(x)= \cdots =\varepsilon_{i+\lceil n^{\alpha+\varepsilon}\rceil}(x)= 0 \ \ \text{for infinitely many} \ n.
\]
Put $a_0=1/2$ and $a_{n+1}=T^{-1}(a_n) \cap [0, 1/2)$. If $\varepsilon_{i}(x)=0$ for all  $1\leq i \leq \lceil n^{\alpha+\varepsilon}\rceil$, then it follows from Lemma 3.2 of \cite{Lsv99} that
\[
x\leq a_{\lceil n^{\alpha+\varepsilon}\rceil} \leq C\frac{1}{\lceil n^{1+\frac{\varepsilon}{\alpha}}\rceil}.
\]
Therefore, we obtain
\begin{equation}\label{mu1}
\mu\{x: \ \varepsilon_{i}(x)=0 \ \text{for all}  \ 0\leq i \leq \lceil n^{\alpha+\varepsilon}\rceil \ \ i.o. \}=0.
\end{equation}

For the second case, it means that there exist $1\leq i \leq n-\lceil n^{\alpha+\varepsilon}\rceil$ and infinitely many $n$ such that
 \[
 T^{i-1}(x)\in I_1 \cap T^{-1}[0, a_{\lceil n^{\alpha+\varepsilon}\rceil}] \subset I_1 \cap T^{-1}[0, a_{\lceil i^{\alpha+\varepsilon}\rceil}].
 \]
Note that $n^{\alpha+\varepsilon}$ is increasing to infinity as $n$ goes to infinity. So, we obtain that either $\varepsilon_{i}(x)=1$ and $\varepsilon_j(x)=0$ for all $j>i$, or $T^{i-1}(x) \in I_1 \cap T^{-1}[0, a_{\lceil i^{\alpha+\varepsilon}\rceil}]$ for infinitely many $i$. It is easy to show that
\begin{equation}\label{mu2}
\mu\{x: \varepsilon_{i}(x)=1 \ \text{and} \ \varepsilon_{j}(x)=0 \ \text{for all} \ j\geq i \}=0.
\end{equation}
Since
\begin{align*}
&\sum_{i=1}^{\infty} \mu\{x:\  \ T^{i-1}(x) \in I_1 \cap T^{-1}[0, a_{\lceil i^{\alpha+\varepsilon}\rceil}] \ \} \\
&=\sum_{i=1}^{\infty} \mu\{I_1 \cap T^{-1}[0, a_{\lceil i^{\alpha+\varepsilon}\rceil}]\}\leq \sum_{i=1}^{\infty} \frac{\tilde{C}}{2i^{1+\frac{\varepsilon}{\alpha}}} < \infty,
\end{align*}
 by Borel-Cantelli lemma, it follows that
\begin{equation}\label{mu3}
\mu \{x: T^{i-1}(x) \in I_1 \cap T^{-1}[0, a_{\lceil i^{\alpha+\varepsilon}\rceil}] \  \text{ for infinitely many} \ i \}=0.
\end{equation}
Combining (\ref{mu1}, \ref{mu2}, \ref{mu3}), we obtain
\[
\mu\{x\in [0,1): r_n^{\alpha}(x)\geq n^{\alpha+\varepsilon} \ i.o. \}=0.
 \]
 Therefore, for $\mu$-almost all $x\in [0,1)$,
 \[
 \limsup_{n\to \infty} \frac{\log r_n^{\alpha}(x)}{\log n} \leq \alpha+\varepsilon.
 \]
 Letting $\varepsilon \to 0$, we complete the proof of Part I.

{\bf Part II.} We will prove that
for $\mu$-almost all $x \in [0,1)$,
\[
\liminf_{n \to \infty} \frac{\log r^\alpha_n(x)}{\log n} \geq \alpha.
\]

We  shall use the following statement. Let $\xi=\{I_0, I_1 \}$ be a partition of $[0, 1)$ and $\xi_n= \xi \bigvee T^{-1}\xi\bigvee \dots \bigvee T^{-n+1}\xi$, where $\xi \bigvee \eta=\{A \cap B: A\in \xi, B\in \eta\}$. In \cite{lesHu04}, Hu showed that there exists $C_1>0$ and $l>0$ such that for any $m\geq 0$ and $A\in\xi_m$ and for any measurable set $B\subset [0, 1) $,
\begin{equation}\label{equ:polynimal}
|\mu(A \cap T^{-n-m}B) - \mu(A) \mu(B)| \leq \frac{C_1 m^{1/\alpha-1}}{(n-l)^{1/\alpha -1}} \mu(A) \mu(B), \  \text{for}\ n>l.
\end{equation}
Then by the invariance of $\mu$, we obtain
\[
|\mu(A^{c} \cap T^{-n-m}B) - \mu(A^c) \mu(B)| \leq \frac{C_1 m^{1/\alpha-1}}{(n-l)^{1/\alpha -1}} \mu(A) \mu(B),
\]
where $A^{c}:=[0, 1) \setminus A$. It implies that
\begin{equation}\label{hu}
\mu(A^{c} \cap T^{-n-m}B) \leq \left( \frac{C_1 m^{1/\alpha-1}}{(n-l)^{1/\alpha -1}} \mu(A) + \mu(A^c)   \right) \mu(B).
\end{equation}

For any $x \in [0,1)$ and $m,n \in \mathbb{N}$ with $m<n$, we define
\[
r^\alpha_{m,n}(x):= \max\{k \geq 0: \varepsilon_{i+1}(x)=\cdots= \varepsilon_{i+k}(x)=0\ \text{for some}\ m-1\leq i \leq n-k\}.
\]
Thus, $r^\alpha_{1,n} = r^\alpha_{n}$.
For any $0<\varepsilon<\alpha$, write $t_n:= \lceil n^{\alpha -\varepsilon}\rceil$ and $k_n: = \lfloor \frac{n}{t^{1+\varepsilon}_n}\rfloor$.
Let
\[
E = \{x\in [0,1): r^\alpha_{t_n}(x) < t_n\}
\]
and
\[
F_j =\{x \in [0,1): r^\alpha_{jt^{1+\varepsilon}_n+1, jt^{1+\varepsilon}_n+t_n}(x) < t_n\},\ \ 1 \leq j \leq k_n-1.
\]
Then
\begin{align*}
\{x \in [0,1): r^\alpha_n(x) < t_n\} &\subseteq E \bigcap F_1 \bigcap \cdots \bigcap F_{k_n-1}\\
&= E \bigcap T^{-t^{1+\varepsilon}_n}(E\cap F_1\cap \cdots \cap F_{k_n-2}).
\end{align*}

Let \[
I_{k}(0,\cdots,0):=\{x \in [0,1):\varepsilon_1(x) = \cdots = \varepsilon_k(x)=0\}\ \text{for some}\ k \geq 1.
\]
Then $I_{t_n}(0,\cdots,0) \in \xi_{t_n}$ and  ${I_{t_n}(0,\cdots,0)}= E^c$. By (\ref{mu}), the similar arguments in Section 6.2 of \cite{lesYoung} imply that
\[
\widetilde{C}^{-1} \frac{1}{t^{1/\alpha -1}_n}\leq \mu(I_{t_n}(0,\cdots,0)) \leq \widetilde{C}\frac{1}{t^{1/\alpha -1}_n},
\]
where $\widetilde{C}>1$ is a constant. For $n$ big enough so that $t_n^{1+\varepsilon} -t_n > l$,
combining this with (\ref{hu}), we deduce that
\begin{align}\label{1}
\mu\{x \in [0,1): r^\alpha_n(x) < t_n\} &\leq \mu(E \bigcap T^{-t^{1+\varepsilon}_n}(E\cap F_1\cap \cdots \cap F_{k_n-2})) \notag\\
&\leq \left(  \frac{C_1t_n^{1/\alpha-1 }}{(t_n^{1+\varepsilon}-t_n -l)^{1/\alpha-1}}\mu(I_{t_n}(0,\dots, 0))+ \mu(E) \right) \notag\\
&~~~~\times \mu (E\cap F_1\cap \cdots \cap F_{k_n-2}) \notag\\
& \leq \cdots \cdots \notag\\
&\leq \left(  \frac{C_1 \widetilde{C} }{(t_n^{1+\varepsilon}-t_n -l)^{1/\alpha-1}}  + 1-  \frac{\widetilde{C}^{-1} }{t^{1/\alpha -1}_n} \right)^{k_n-1}.
\end{align}

Using the inequality: $e^{-x} \geq 1-x$ for any $x \geq 0$, we have
\begin{align*}
\mu\{x \in [0,1): r^\alpha_n(x) < t_n\} & \leq \left(1- \frac{\widetilde{C}^{-1} }{t^{1/\alpha -1}_n} + o(\frac{1 }{t^{1/\alpha -1}_n}) \right)^{k_n-1}\\
&\leq  \exp\left( -\frac{\widetilde{C}^{-1} }{t^{1/\alpha -1}_n} \frac{n}{t^{1+\varepsilon}_n} \right)\\
&=\exp\left(-\widetilde{C}^{-1} n^{-(\frac{\varepsilon }{\alpha}+\varepsilon^2-\alpha\varepsilon)} \right).
\end{align*}

Since $\varepsilon/\alpha + \varepsilon^2 > \alpha \varepsilon$, we eventually get that
\[
\sum_{n \geq 1} \mu\{x \in [0,1): r^\alpha_n(x) < t_n\} <\infty.
\]
By Borel-Cantelli lemma, we conclude that
for $\mu$-almost all $x \in [0,1)$,
\[
\liminf_{n \to \infty} \frac{\log r^\alpha_n(x)}{\log n} \geq \alpha-\varepsilon.
\]
Letting $\varepsilon \to 0^+$, the proof of Part II is completed.

{\bf Part III.} We are concerned in this part with the run length function $R^\alpha_n(x)$.
This part can be proved by the similar way as shown before.
We first show that for $\mu$-almost all $x \in [0,1)$,
\[
\limsup_{n \to \infty} \frac{R^\alpha_n(x)}{\log_2n} \leq 1.
\]
In fact, let $\varepsilon>0$, if $R_n^{\alpha}(x) \geq (1+\varepsilon)\log_2n$ holds for some $x\in [0, 1)$ and $n\in \mathbb{N}$, then there exists $0\leq i \leq n-k$ (where $k=\lceil (1+\varepsilon) \log_2n \rceil$) such that
$$
\varepsilon_{i+1}(x)=\dots =\varepsilon_{i+k}(x)=1.
$$
That is, $T^i(x) \in I_1$, $T^{i+1}(x) \in I_1$, $\dots$, $T^{i+k-1}(x)\in I_1$. Hence
\[
1-\frac{1}{(i+1)^{1+\varepsilon}} \leq 1-\frac{1}{2^k} \leq T^i(x) < 1.
\]
Since $(1+\varepsilon) \log_2n \to \infty$ as $n \to \infty$, $R_n^{\alpha}(x) \geq (1+\varepsilon)\log_2n$ for infinite many $n$ implies that
either $T^i(x)\in [1-\frac{1}{{(i+1)}^{1+\varepsilon}}, 1]$ for infinite many $i$ or there exists an $i\in \mathbb{N}$ such that $\varepsilon_j=1$ for all $j\geq i$. Then it follows from Borel-Cantelli lemma that
\[
\mu\{x\in [0,1): R_n^{\alpha}(x) \geq (1+\varepsilon)\log_2n \ i.o.\ \}=0.
\]
Therefore, for $\mu$-almost all $x\in [0,1)$,
\[
\limsup_{n \to \infty} \frac{R_n^{\alpha}(x)}{\log_2n}\leq 1+\varepsilon.
\]
Taking $\varepsilon \to 0$, we obtain the proof of the first part.

Next we will show for $\mu$-almost all $x \in [0,1)$,
\[
\liminf_{n \to \infty} \frac{R^\alpha_n(x)}{\log_2n} \geq 1.
\]
We  shall use the following facts. For any interval $A \subseteq (1/2,1)$ and interval $B \subseteq (1/2,1)$,
\begin{equation*}
|\mu(A\cap T^{-n}B) - c_n\mu(A)\mu(B)| \leq \frac{1}{n^{1/\alpha}}\mu(B),\ \ \forall n\geq 1,
\end{equation*}
where $c_n =  1+ \frac{s}{n^{1/\alpha-1}} + o(\frac{1}{n^{1/\alpha-1}})$ for some non-zero constant $s$, see Eq. (1.3) of Gou\"{e}zel \cite{lesGou07}. Then
\begin{equation}\label{mixing}
\mu(A^c\cap T^{-n}B) \leq \left(c_n\mu(A^c)+ \frac{1}{n^{1/\alpha}}\right)\mu(B),\ \ \forall n\geq 1.
\end{equation}

For any $x \in [0,1)$ and $m,n \in \mathbb{N}$ with $m<n$, we define
\[
R^\alpha_{m,n}(x):= \max\{k \geq 0: \varepsilon_{i+1}(x)=\cdots= \varepsilon_{i+k}(x)=1\ \text{for some}\ m-1\leq i \leq n-k\}.
\]
Thus, $R^\alpha_{1,n} = R^\alpha_{n}$.
For $\varepsilon>0$ small such that $\frac{\alpha}{1-\alpha} - \varepsilon -\varepsilon^2 \frac{1}{1-\alpha}>0$, write $t_n:= \lceil (1-\varepsilon) \log_2n \rceil$, and $l_n= \lfloor  2^{\frac{\alpha t_n (1+\varepsilon)}{1-\alpha}} \rfloor$. Set $k_n: = \lfloor \frac{n}{l_n}\rfloor$.
Let
\[
E = \{x\in [0,1): R^\alpha_{t_n}(x) < t_n\}
\]
and
\[
F_j =\{x \in [0,1): R^\alpha_{jl_n+1, jl_n+t_n}(x) < t_n\},\ \ 1 \leq j \leq k_n-1.
\]
Then
\begin{align*}
\{x \in [0,1): R^\alpha_n(x) < t_n\} &\subseteq E \bigcap F_1 \bigcap \cdots \bigcap F_{k_n-1}\\
&= E \bigcap T^{-l_n}(E\cap F_1\cap \cdots \cap F_{k_n-2}).
\end{align*}

Note that $I_{t_n}(1,\cdots,1) \subseteq [1/2, 1)$ and $I_{t_n}(1,\cdots,1)=E^c$, where
\[
I_{t_n}(1,\cdots,1):=\{x \in [0,1):\varepsilon_1(x) = \cdots = \varepsilon_{t_n}(x)=1\}.
\]
By (\ref{mu0}), it follows that
\[
\widetilde{C}^{-1} \frac{1}{2^{t_n}}\leq \mu(I_{t_n}(1,\cdots,1)) \leq \widetilde{C}\frac{1}{2^{t_n}},
\]
where $\widetilde{C}>1$ is a constant.
Combing this with (\ref{mixing}), we deduce that
\begin{align*}
\mu\{x \in [0,1): R^\alpha_n(x) < t_n\} &\leq \mu(E \bigcap F_1 \bigcap \cdots \bigcap F_{k_n-1}) \notag\\
&\leq \left(c_{l_n}\mu(E)+ \frac{1}{l_n^{1/\alpha}} \right) \mu (E\cap F_1\cap \cdots \cap F_{k_n-2}) \notag\\
& \leq \cdots \cdots \notag\\
&\leq \left(c_{l_n}\mu(E)+ \frac{1}{l_n^{1/\alpha}} \right)^{k_n-1} \\
 &\leq \left(\left(1+ \frac{s}{2^{t_n(1+\varepsilon)}} + o(\frac{1}{2^{t_n(1+\varepsilon)}})\right)\left(1- \widetilde{C}^{-1} \frac{1}{2^{t_n}}\right)+ \frac{1}{2^{\frac{t_n(1+\varepsilon)}{1-\alpha}}}\right)^{k_n-1}\\
 & \leq \left(1- \frac{\widetilde{C}^{-1} }{2^{t_n}} + o(\frac{1}{2^{t_n}}) \right)^{k_n-1} \\
 & \leq \exp\left\{ -\frac{\widetilde{C}^{-1}}{ 2^{t_n}} \frac{n}{2^{\frac{t_n \alpha(1+\varepsilon)}{1 -\alpha}}}    \right\}.
\end{align*}
Then the choosing of $\varepsilon$ implies that

\[
\sum_{n \geq 1} \mu\{x \in [0,1): R^\alpha_n(x) < t_n\} <\infty,
\]
by Borel-Cantelli lemma, we conclude that
for $\mu$-almost all $x \in [0,1)$,
\[
\liminf_{n \to \infty} \frac{R^\alpha_n(x)}{\log_2 n} \geq 1-\varepsilon.
\]
Letting $\varepsilon \to 0^+$, we complete the proof.

\end{proof}

\section{Discussions for further research}
In this paper we have studied the run length function for the intermittent maps $T_{\alpha}$ with $0<\alpha<1$. A natural question is to estimate the run length function for intermittent maps with $\alpha\geq1$. Notice that an intermittent map $T_{\alpha}$ admits a unique ergodic invariant probability measure absolutely continuous with respect to Lebesgue measure if and only if $0<\alpha<1$, and in this case, the decay of correlation is polynomial. As we have mentioned in the Introduction, such polynomial decay of correlations and the classical Borel-Cantelli lemma are two fundamental ingredients in our proofs of the main results.


The difficulties in obtaining pointwise estimation of the run length function for intermittent maps with $\alpha \geq 1$ are twofold. Firstly, the decay of correlation of such maps do not hold. To be more precise, when $\alpha\geq 1$, each intermittent map admits a unique (up to scaling) $\sigma$-finite (but not finite), absolutely continuous invariant measure $\mu$. In these infinite measure settings, although one may still obtain certain mixing rate $\int_{[0,1]}fg\circ T^n_\alpha d\mu$ for some reasonably well-behaved observables $f,g$. (See e.g.,\cite{MT12, Ter16} for the concrete statements), there is no analogous result as in \eqref{equ:polynimal} to the best of our knowledge. Secondly, the classical Borel-Cantelli lemma does not hold in the infinite measure setting. Therefore, our methods in this paper are inapplicable to obtain the corresponding results for $\alpha\geq 1$ directly, and we expect to develop more new methods and techniques to overcome these difficulties.

\section*{Acknowledgments}
We would like to thank Prof. Manfred Denker and Prof. Huyi Hu for useful discussions, particularly for the proof of part I of our main theorem. This work is partially supported by grants from National Natural Science Foundation of China (Nos. 11701200, 11671395, 11801591, 11871262), Hubei Chenguang Talented Youth Development Foundation 2017 (Nos. 0106011025), and Ky and Yu-Fen Fan traveling award, National Science Foundation 2018.

\end{document}